\documentclass{amsart}

\usepackage{graphicx}
\usepackage{pinlabel}

\newtheorem{theorem}{Theorem}
\newtheorem{lemma}[theorem]{Lemma}
\newtheorem{corollary}[theorem]{Corollary}

\begin{document}

\title{Local rigidity of inversive distance circle packing}

\author{Ren Guo}

\address{School of Mathematics, University of Minnesota, Minneapolis, MN, 55455}

\email{guoxx170@math.umn.edu}

\thanks{}



\begin{abstract}
A Euclidean (or hyperbolic) circle packing on a closed triangulated surface with prescribed inversive distance is locally determined by its cone angles. We prove this by applying a variational principle.
\end{abstract}

\maketitle

\section{Introduction}

\subsection{Andreev-Thurston Theorem}

In his work on constructing hyperbolic metrics on 3-manifolds, W. Thurston (\cite{th}, Chapter 13) studied a Euclidean (or hyperbolic) circle packing on a closed triangulated surface with prescribed intersection angles.
Thurston's work generalizes Andreev's result of circle packing on a sphere \cite{a1,a2}. The special case of tangent circle packing (every intersection angle is zero) on sphere was obtained by Koebe \cite{k}. 

Suppose $(\Sigma,T)$ is a closed triangulated surface so that
$V,E,F$ are sets of all vertices, edges and triangles in $T.$ We identify vertexes of $T$ with indexes, edges of $T$ with pairs of indexes and triangles of $T$ with triples of indexes. This means $V=\{1,2,...|V|\}, E=\{ij\ |\ i,j\in V\}$ and $F=\{\triangle ijk\ |\ i,j,k\in V\}.$ Fix a vector $\Theta\in \mathbb{R}^{|E|}$ indexed by the set of edges $E$, such that $\Theta_{ij}\in [0,\frac\pi2]$ for the each $ij\in E.$ This vector is call a \it weight \rm on $(\Sigma,T)$.

A circle packing on a closed triangulated weighted surface $(\Sigma,T,\Theta)$ is a configuration $\{c_i,i\in V\}$ of circles such that the intersection angle of $c_i$ and $c_j$ is $\Theta_{ij}$ which is the prescribed weight on the edge $ij.$ To obtain a Euclidean (or hyperbolic) circle packing for $(\Sigma,T,\Theta),$ we start with a \it radius vector \rm $r=(r_1,r_2,...,r_{|V|})\in \mathbb{R}^{|V|}_{>0}$ which assigns each vertex $i\in V$ a positive number $r_i.$ A radius vector $r$ produces a Euclidean (or hyperbolic) cone metric on the surface as follows.

\begin{figure}[ht!]
\labellist\small\hair 2pt
\pinlabel $i$ at 177 278
\pinlabel $j$ at 86 126
\pinlabel $k$ at 305 126
\pinlabel $p$ at 102 212
\pinlabel $r_i$ at 132 252
\pinlabel $r_j$ at 91 172
\pinlabel $\Theta_{ij}$ at 82 234
\pinlabel $\Theta_{jk}$ at 178 50
\pinlabel $\Theta_{ki}$ at 273 269
\pinlabel $l_{ij}$ at 141 182
\pinlabel $l_{jk}$ at 205 139
\pinlabel $l_{ki}$ at 229 182

\endlabellist
\centering
\includegraphics[scale=0.5]{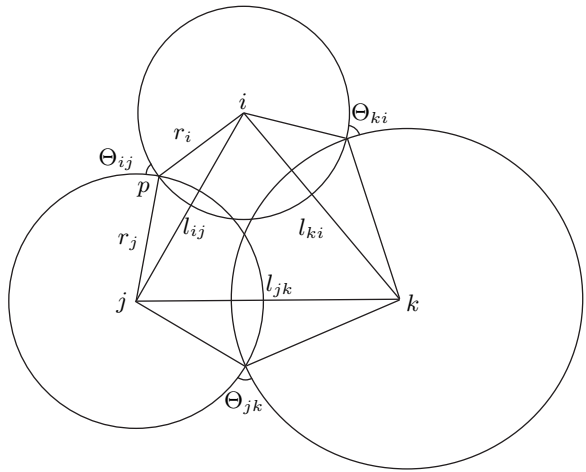}
\caption{Andreev-Thurston circle packing}
\label{fig:overlap}
\end{figure}

In Figure \ref{fig:overlap}, consider a topological triangle $\triangle ijk.$
One can construct a Euclidean (or hyperbolic)
triangle $\triangle ijp$ such that the edges $ip,jp$
have lengths $r_i,r_j$ respectively and the angle at $p$
is $\pi-\Theta_{ij}$. Let $l_{ij}$ be the length of the edge $ij$
in the Euclidean (or hyperbolic) triangle $\triangle ijp$ which is a
function of $r_i,r_j.$ In fact, in Euclidean geometry
$$l_{ij}=\sqrt{r^2_i+r^2_j+2r_ir_j\cos \Theta_{ij}}$$ or in hyperbolic geometry
$$l_{ij}=\cosh^{-1}(\cosh r_i\cosh r_j+\cos \Theta_{ij}\sinh r_i\sinh r_j).$$
Similarly, we obtain $l_{jk},l_{ki}.$

As observed by Thurston \cite{th}, under the assumption that $\Theta_{ij}\in [0,\frac\pi2]$ for each $ij\in E$, $l_{ij},l_{jk}$ and $l_{ki}$ can be realized as edge lengths of a Euclidean (or hyperbolic) triangle $\triangle
ijk$ in $F$. Gluing all of these Euclidean (or hyperbolic) triangles in $F$ produces a Euclidean (or hyperbolic) cone metric on $(\Sigma,T)$ with possible cones at vertexes of $T$. On this surface with a cone metric, drawing a circle centered at vertex $i$ with radius $r_i$ for each vertex $i\in V,$ we obtain a circle packing with prescribed intersection angles.

At each vertex $i\in V$, the cone angle $a_i$ is the sum of all inner angles having vertex $i.$ Thus cone angles $(a_1,a_2,...,a_{|V|})$ are functions of the radius vector $r=(r_1,r_2,...,r_{|V|})$. We try to find a radius vector when cone angles are given. Andreev-Thurston Theorem answer the question about the existence and uniqueness of solutions of radius vector when cone angles are given.

\begin{theorem}[Andreev-Thurston] For any closed triangulated weighted surface
$(\Sigma,T,\Theta)$ with $\Theta_{ij}\in [0,\frac\pi2]$ for each $ij\in E$, a Euclidean (or hyperbolic) circle
packing for $(\Sigma,T,\Theta)$ is determined by its cone angles up to Euclidean similarity (or hyperbolic isometry).
Furthermore, the set of all possible cone angles form an open convex polytope in
$\mathbb{R}^{|V|}.$
\end{theorem}

For a proof, see Thurston \cite{th}, Marden-Rodin \cite{mr}, Colin de Verdi\'ere
\cite{cdv}, He \cite{h}, Chow-Luo \cite{cl}, Stephenson \cite{s}.

\subsection{Inversive distance}

H. S. M. Coxster \cite{c} introduced the notion of inversive distance to describe the relationship between two circles in a M\"obius plane.
The notion of inversive distance generalizes the notion of intersection angle of two circles.

The inversive distance between two circles is independent of geometry, i.e., the inversive distance of two circles is the same under spherical, Euclidean and hyperbolic metrics. This generalizes the fact that the intersection angle of two circles is the same under spherical, Euclidean and hyperbolic metrics.

In the Euclidean plane, consider two circles $c_1, c_2$ with radii $r_1, r_2$ respectively. In the relevant case of this paper, we assume that $c_i$ does not contain $c_j$ for $\{i,j\}=\{1,2\}.$ If the distance between their center is $l$, the inversive distance between $c_1, c_2$ is given by the formula
$$I(c_1,c_2)=\frac{l^2-r_1^2-r_2^2}{2r_1r_2}.$$

In hyperbolic plane, consider two circles $c_1, c_2$ with radii $r_1, r_2$ respectively. We assume that $c_i$ does not contain $c_j$ for $\{i,j\}=\{1,2\}.$ If the distance between their center is $l$, the inversive distance between $c_1, c_2$ is given by the formula
$$I(c_1,c_2)=\frac{\cosh l -\cosh r_1\cosh r_2}{\sinh r_1\sinh r_2}.$$

\begin{figure}[ht!]
\labellist\small\hair 2pt
\pinlabel $\mathbb{H}^3$ at 271 215
\pinlabel $\mathbb{C}$ at 332 169
\pinlabel $\mathbb{H}^2$ at 82 152
\pinlabel $c_1$ at 109 127
\pinlabel $c_2$ at 245 116
\pinlabel $\bar{c}_1$ at 132 77
\pinlabel $\bar{c}_2$ at 222 73
\pinlabel $(0,0,-1)$ at 220 2

\endlabellist
\centering
\includegraphics[scale=0.7]{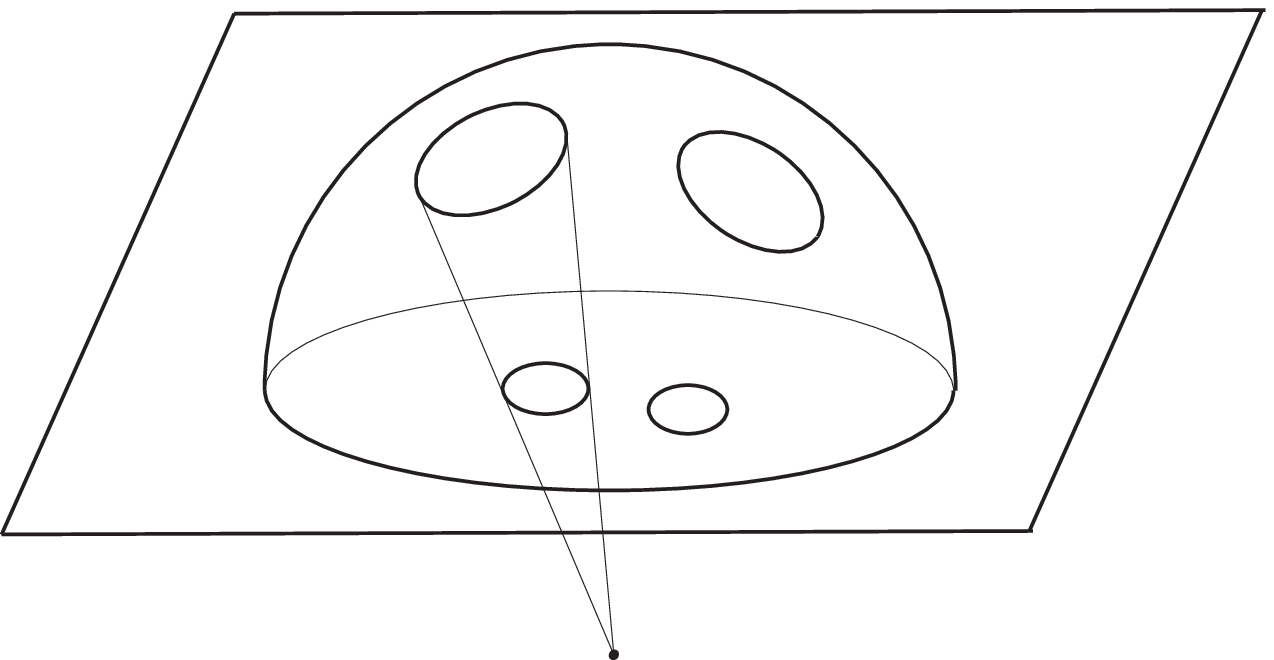}
\caption{the stereographic projection}
\label{fig:projection}
\end{figure}

These two formulas are related by a stereographic projection. In Figure \ref{fig:projection}, consider the upper half space model of the hyperbolic space $\mathbb{H}^3$. The hyperbolic plane $\mathbb{H}^2$ is realized as a totally geodesic plane in $\mathbb{H}^3$, where it is a unit hemisphere perpendicular to the infinite boundary $\mathbb{C}$. The stereographic projection is the projection staring from point $(0,0,-1),$ which sends objects in $\mathbb{H}^2$ into $\mathbb{C}$.
Let $c_1,c_2$ be two circles in $\mathbb{H}^2$ with radii $r_1,r_2$ and the distance between centers $l$. The stereographic projection produces two circles $\bar{c}_1,\bar{c}_2$ on the Euclidean plane $\mathbb{C}$ with radii $\bar{r}_1,\bar{r}_2$ and the distance between centers $\bar{l}$. Then we have
$$I(c_1,c_2)=\frac{\cosh l -\cosh r_1\cosh r_2}{\sinh r_1\sinh r_2}=\frac{\bar{l}^2-\bar{r}_1^2-\bar{r}_2^2}{2\bar{r}_1\bar{r}_2}=I(\bar{c}_1,\bar{c}_2).$$
We say that the stereographic projection preserves the inversive distance.

When $I(c_1,c_2)\in (-1,0),$ the circles $c_1, c_2$ intersect with an intersection angle $\arccos I(c_1,c_2)\in (\frac{\pi}2,\pi).$ When $I(c_1,c_2)\in [0,1),$ the circles $c_1, c_2$ intersect with an intersection angle $\arccos I(c_1,c_2)\in (0,\frac\pi2].$ When $I(c_1,c_2)=1,$ the circles $c_1, c_2$ are tangent to each other. When $I(c_1,c_2)\in (1,\infty),$ the circles $c_1, c_2$ are separated.

\subsection{Inversive distance circle packing}

Motivated by the application to discrete conformal mappings, Bowers-Stephenson \cite{bs} introduced inversive distance circle packing which generalizes Andreev-Thurston's intersection angle circle packing. See Stephenson \cite{s} (page 331) and Bowers-Hurdal \cite{bh} for more information.

We will use the same notation as in Andreev-Thurston's intersection angle circle packing. For a closed triangulated surface $(\Sigma,T)$, fix a weight $I\in \mathbb{R}^{|E|}$ with $I_{ij}\in [0,\infty)$ for each $ij\in E.$ In this paper, we will not consider the case of obtuse intersection angles. (Rivin \cite{r} and Leibon \cite{le} obtained results for certain circle patterns on surfaces allowing obtuse intersection angles.)
An inversive distance circle packing for $(\Sigma,T,I)$ is a configuration of circles $\{c_i,i\in V\}$ such that the inversive distance of two circles $c_i,c_j$ is the prescribed number $I_{ij}$ for the edge $ij.$
To obtain a Euclidean (or hyperbolic) circle packing for $(\Sigma,T,I),$ we start with a radius vector $r=(r_1,r_2,...,r_{|V|})\in \mathbb{R}^{|V|}_{>0}$ which assigns each vertex $i\in V$ a positive number $r_i.$ A radius vector $r$ produces a Euclidean (or hyperbolic) cone metric on the surface as follows.

Consider a topological triangle
$\triangle ijk\in F$. Assign a number $l_{ij}$ to the edge $ij\in E,$ where $l_{ij}$ is the distance of two circles with radii $r_i,r_j$ and inversive distance $I_{ij}.$ Therefore $l_{ij}$ is a function of $r_i, r_j.$ In fact, in Euclidean geometry
$$l_{ij}=\sqrt{r^2_i+r^2_j+2r_ir_jI_{ij}}$$ or in hyperbolic geometry
$$l_{ij}=\cosh^{-1}(\cosh r_i\cosh r_j+I_{ij}\sinh r_i\sinh r_j).$$
Similarly, we obtain $l_{jk},l_{ki}.$

To guarantee that $l_{ij},l_{jk},l_{ki}$ satisfy the triangle inequality for each triangle $\triangle ijk \in F,$ a vector $r=(r_1,r_2,...,r_{|V|})$ must satisfy a certain conditions. For details, see section 2.1 and 3.1.

For a radius vector $r$, if $l_{ij},l_{jk},l_{ki}$ can be realized as edge lengths of a Euclidean (or hyperbolic) triangle for each $\triangle ijk \in F,$ gluing all of these Euclidean (or hyperbolic) triangles in $F$ produces a Euclidean (or hyperbolic) cone metric on $(\Sigma,T)$ with possible cones at vertexes of $T$. On this surface with a cone metric, drawing a circle centered at vertex $i$ with radius $r_i$ for each vertex $i\in V,$ we obtain a circle packing with prescribed inversive distance.

The cone angles $(a_1,a_2,...,a_{|V|})$ are functions of the radius vector $r=(r_1,r_2,...,r_{|V|}).$ We are interested in the solutions of radius vector when cone angles are given.

Our main theorem is the local rigidity of inversive distance circle packing.

\begin{theorem} \label{thm:main} For any closed triangulated weighted surface
$(\Sigma,T, I)$ with $I_{ij}\in [0,\infty)$ for each $ij\in E,$ a Euclidean (or hyperbolic) inversive distance circle packing for $(\Sigma,T,I)$ is locally determined by its cone angle up to Euclidean similarity (or hyperbolic isometry).
\end{theorem}

The local rigidity of Andreev-Thurston Theorem is generalized to the case of inversive distance circle packing. In our case, we are not able to obtain either the global rigidity or the description of the space of cone angles as in Andreev-Thurston Theorem. For discussion, see section 2.1 and 3.1.

In section 2, we prove Theorem \ref{thm:main} in Euclidean geometry. In section 3, we prove Theorem \ref{thm:main} in hyperbolic geometry.

\subsection{Variational principle}

Theorem \ref{thm:main} is proved by applying a variational principle. The energy function of this variational principle is derived from the derivative of the cosine law of a Euclidean (or hyperbolic) triangle. Andreev and Thurston obtained their results without using variational principles. The variational approach to circle packing was first introduced by Colin de Verdi\'ere
\cite{cdv}. Since then, many works on variational principles on
circle packing or circle pattern have appeared. For example, see Br\"agger \cite{b}, Rivin \cite{r}, Leibon \cite{le}, Chow-Luo \cite{cl}, Bobenko-Springborn \cite{b-s}, Guo \cite{g}, Springborn \cite{sp}, Luo \cite{l}, Guo-Luo \cite{gl} and others.

Applying a variational principle, Colin de Verdi\'ere
\cite{cdv} proved Andreev-Thurston Theorem in the special case of tangent circle packing, i.e., $\Theta_{ij}=0$ for each $ij\in E$. Chow-Luo \cite{cl} proved Andreev-Thurston Theorem of intersection angle circle packing, i.e., $\Theta_{ij}\in [0,\frac\pi2]$ for each $ij\in E$.
Variational principles for polyhedral surfaces including the topic of circle packing are studied systematically in Luo \cite{l}. Many energy functions are derived from the cosine law and its derivative. Tangent circle packing is generalized to tangent circle packing with a family of discrete curvature. For exposition of this work, see also Dai-Gu-Luo \cite{dgl}. We follow this variational approach to prove Theorem \ref{thm:main} of inversive distance circle packing, i.e., $I_{ij}\in [0,\infty)$ for each $ij\in E$.

\section{Circle packing in Euclidean geometry}

\subsection{Space of radius vectors}

As stated in Introduction, to obtain a Euclidean cone metric, we need conditions on a radius vector. In this section we investigate the space of radius vectors for one triangle in Euclidean geometry.
Consider a topological triangle
$\triangle ijk\in F$. Assign a number $l_{ij}$ to the edge $ij\in E,$ where $l_{ij}$ is the distance of the centers of two circles with radii $r_i,r_j$ and inversive distance $I_{ij}.$ Therefore $l_{ij}$ is a function of $r_i, r_j.$ In fact, in Euclidean geometry
$$l_{ij}=\sqrt{r^2_i+r^2_j+2r_ir_jI_{ij}}.$$ Similarly, we obtain $l_{jk},l_{ki}.$ 

\begin{lemma}\label{thm:e-r-space} The set $\mathcal{R}_E^{ijk}:=\{(r_i, r_j,r_k)\in \mathbb{R}^3_{>0}\ |\ l_{ij},l_{jk},l_{ki}$ satisfy triangle inequality\} is a connected and simply connected open subset in $\mathbb{R}^3_{>0}.$
\end{lemma}

\begin{proof} Let $\rho_E:\mathbb{R}^3_{> 0}\to \mathbb{R}^3_{> 0}$ be the map sending $(r_i,r_j,r_k)$ to $(l_{ij},l_{jk},l_{ki}).$ We claim that $\rho_E$ is a smooth embedding. First the determinant of the Jacobian matrix of $\rho_E$ is nonzero. Next, we will show $\rho_E$ is one-to-one.

That means we need to show if the following equations has a solution $(r_i,r_j,r_k)$, then the solution is unique.
$$r_j^2+r_k^2+2I_{jk}r_jr_k=l_{jk}^2,$$
$$r_k^2+r_i^2+2I_{ki}r_kr_i=l_{ki}^2,$$
$$r_i^2+r_j^2+2I_{ij}r_ir_j=l_{ij}^2.$$

Solving the first equation for $r_j$, we get $r_j=-I_{jk}r_k\pm\sqrt{(I_{jk}^2-1)r_k^2+l_{jk}^2}.$ Since $I_{jk}>0$ and $r_j>0,$ we must have $r_j=-I_{jk}r_k+\sqrt{(I_{jk}^2-1)r_k^2+l_{jk}^2}.$ Hence $r_j$ is uniquely determined by $r_k.$ By the same argument $r_i=-I_{ki}r_k+\sqrt{(I_{ki}^2-1)r_k^2+l_{ki}^2}.$

After substituting $r_i,r_j$ in the third equation, we need to show it has a unique solution for $r_k$.
Consider the function $$f(t)=r_i^2(t)+r_j^2(t)+2I_{ij}r_i(t)r_j(t)-l_{ij}^2$$ which is obtained by replacing $r_k$ by $t$ for simplicity of notations.
Now $$f'(t)=(2r_i+2I_{ij}r_j)r_i'(t)+(2r_j+2I_{ij}r_i)r_j'(t).$$
And we have $$r_i'(t)=-I_{ki}+\frac{(I_{ki}^2-1)t}{\sqrt{(I_{ki}^2-1)t^2+l_{ki}^2}}.$$

If $I_{ki}\in[0,1],$ then each term in $r_i'(t)$ is non-positive. Thus $r_i'(t)<0.$

If $I_{ki}\in(1,\infty),$ then
$$r_i'(t)<-I_{ki}+\frac{(I_{ki}^2-1)t}{\sqrt{(I_{ki}^2-1)t^2}}=-I_{ki}+\sqrt{I_{ki}^2-1}<0.$$

By the same argument, we have $r_j'(t)<0.$ Since $2r_i+2I_{ij}r_j>0$ and $2r_j+2I_{ij}r_i>0,$ then $f'(t)<0$ holds. Therefore $f(t)=0$ has a unique solution.

To determine the image $\rho_E(\mathbb{R}^3_{> 0}),$ we consider the image of the boundary of $\mathbb{R}^3_{\geq 0}.$ If $r_i=0,$ then $r_k=l_{ki}, r_j=l_{ij}.$ Thus $$l_{jk}^2=r_j^2+r_k^2+2I_{jk}r_jr_k=l_{ij}^2+l_{ki}^2+2I_{jk}l_{ij}l_{ki}.$$ On the other hand, if $r_i>0,$ then $r_k<l_{ki}, r_j<l_{ij}.$ We have $$l_{jk}^2<l_{ij}^2+l_{ki}^2+2I_{jk}l_{ij}l_{ki}.$$ Therefore the image $\rho_E(\mathbb{R}^3_{> 0})$ is the subset of vectors $(l_{ij},l_{jk},l_{ki})$ satisfying
$$l_{jk}^2< l_{ij}^2+l_{ki}^2+2I_{jk}l_{ij}l_{ki},$$
$$l_{ki}^2< l_{jk}^2+l_{ij}^2+2I_{ki}l_{jk}l_{ij},$$
$$l_{ij}^2< l_{ki}^2+l_{jk}^2+2I_{ij}l_{ki}l_{jk}.$$
In fact, $\rho_E(\mathbb{R}^3_{> 0})$ is a cone bounded by three surfaces. In $\mathbb{R}^3_{\geq 0}$, the intersection of any of the two surfaces is a straight line. The three lines of intersection are 
$$l_{ij}=0 \ \mbox{and}\ l_{jk}=l_{ki},$$
$$l_{jk}=0 \ \mbox{and}\ l_{ki}=l_{ij},$$
$$l_{ki}=0 \ \mbox{and}\ l_{ij}=l_{jk}.$$
 
Let $\mathcal{L}$ be the subset of $\mathbb{R}^3_{> 0}$ formed by vectors $(l_{ij},l_{jk},l_{ki})$ satisfying the triangle inequality. Therefore $\mathcal{L}$ is a cone bounded by three planes. The three lines of intersection of the three planes are the same as that of the boundary surfaces of $\rho_E(\mathbb{R}^3_{> 0}).$

We claim that $\rho_E(\mathbb{R}^3_{> 0})\cap \mathcal{L}$ is a cone. Therefore it is connected and simply connected.

In fact, if $I_{ij}\in [0,1],$ then, in $\mathbb{R}^3_{> 0}$, we have
$$\{(l_{ij},l_{jk},l_{ki})\ |\ l_{ij}^2< l_{ki}^2+l_{jk}^2+2I_{ij}l_{ki}l_{jk}\}\subseteq\{(l_{ij},l_{jk},l_{ki})\ |\ l_{ij}< l_{ki}+l_{jk}\}.$$
If $I_{ij}\in (1,\infty),$ then, in $\mathbb{R}^3_{> 0}$, we have
$$\{(l_{ij},l_{jk},l_{ki})\ |\ l_{ij}^2< l_{ki}^2+l_{jk}^2+2I_{ij}l_{ki}l_{jk}\}\supset\{(l_{ij},l_{jk},l_{ki})\ |\ l_{ij}< l_{ki}+l_{jk}\}.$$

To determine the shape of $\rho_E(\mathbb{R}^3_{> 0})\cap \mathcal{L}$, there are four cases to consider.

If $I_{ij}, l_{jk},l_{ki}\in (1.\infty),$ then $\rho_E(\mathbb{R}^3_{> 0})\cap \mathcal{L}=\mathcal{L}.$

If $I_{ij}\in[0,1], l_{jk},l_{ki}\in (1.\infty),$ then $\rho_E(\mathbb{R}^3_{> 0})\cap \mathcal{L}$ is a cone bounded by one surface and two planes.

If $I_{ij}, l_{jk}\in[0,1],l_{ki}\in (1.\infty),$ then $\rho_E(\mathbb{R}^3_{> 0})\cap \mathcal{L}$ is a cone bounded by two surfaces and one plane.

If $I_{ij}, l_{jk},l_{ki}\in[0,1],$ then $\rho_E(\mathbb{R}^3_{> 0})\cap \mathcal{L}=\rho_E(\mathbb{R}^3_{> 0}).$ 

To finish the proof of the lemma, by definition, $\mathcal{R}_E^{ijk}=\rho_E^{-1}(\rho_E(\mathbb{R}^3_{> 0})\cap \mathcal{L}).$ Since $\rho_E$ is a smooth embedding, $\mathcal{R}_E^{ijk}$ is connected and simply connected.
 
\end{proof}

In the construction of the energy function, we will use a natural variable $u=(u_1,u_2,...,u_{|V|})$ where $u_i=
\ln (r_i).$ Let $\tau_E:\mathbb{R}^3\to \mathbb{R}^3_{>0}$ be the map sending $(u_i,u_j,u_k)$ to $(r_i,r_j,r_k).$ Since $\tau_E$ is a diffeomorphism, we have 

\begin{corollary}\label{thm:e-space}
The set $\mathcal{U}_E^{ijk}:=\tau_E^{-1}(\mathcal{R}_E^{ijk})$ is a connected and simply connected open subset in $\mathbb{R}^3.$
\end{corollary}

To prove Theorem \ref{thm:main}, we construct a concave energy function on $\mathcal{U}_E^{ijk}$ for each triangle $\triangle ijk\in F.$
We obtain the local rigidity of inversive distance circle packing from this energy function. We can not prove the global rigidity due to the fact that $\mathcal{U}_E^{ijk}$ is not always convex in this general setting. In fact, $l_{ij},l_{jk},l_{ki}$ satisfy the triangle inequality is equivalent to
\begin{equation}\label{fml:e-triangle}
(l_{ij}+l_{jk}+l_{ki})(l_{jk}+l_{ki}-l_{ij})(l_{jk}+l_{ij}-l_{ki})(l_{ij}+l_{jk}-l_{ki})>0.
\end{equation}
Substituting the functions $l_{ij},l_{jk},l_{ki}$ in terms of $u_i,u_j,u_k$ into inequality (\ref{fml:e-triangle}), we see that $\mathcal{U}_E^{ijk}$ is defined by the inequality
\begin{multline}\label{fml:e-triangle-u}
\frac{1-I_{jk}^2}{e^{2u_i}}+\frac{1-I_{ki}^2}{e^{2u_j}}+\frac{1-I_{ij}^2}{e^{2u_k}}\\
                              +\frac{2(I_{jk}I_{ki}+I_{ij})}{e^{u_i+u_j}}+\frac{2(I_{ki}I_{ij}+I_{jk})}{e^{u_j+u_k}}+\frac{2(I_{ij}I_{jk}+I_{ki})}{e^{u_k+u_i}}>0.
\end{multline}
For example, when $I_{ij}=I_{jk}=I_{ki}=2,$ the inequality (\ref{fml:e-triangle-u}) determines a subset which is not convex.

In Andreev-Thurston Theorem, the space of cone angle is determined by continuity argument. Since we can not prove the global rigidity, we are not able to determine the space of cone angle in the case of inversive distance circle packing.

\subsection{Energy function}

To prove Theorem \ref{thm:main}, we will construct an energy function on the space of radius vectors. First, we construct an energy function on the space for one triangle. To simplify notations, consider the triangle $\triangle 123\in F$ and the map
$$\begin{array}{ccccccc}
\mathcal{U}_E^{123} & \to & \mathcal{R}_E^{123} & \to &  \mathcal{L} & \to &  \mathbb{R}^3_{>0}\\
(u_1,u_2,u_3) & \mapsto &  (r_1,r_2,r_3) & \mapsto &  (l_1,l_2,l_3) & \mapsto &  (\alpha_1,\alpha_2,\alpha_3)
\end{array}$$
where $r_i=e^{u_i}$ for $i=1,2,3.$ For three fixed positive numbers $I_{12},I_{23},I_{31},$ we have $l_i^2=r_j^2+r_k^2+2I_{jk}r_jr_k$ for $\{i,j,k\}=\{1,2,3\}.$ For simplicity of notations, we set $(l_1,l_2,l_3)=(l_{23},l_{31},l_{12}).$ And $\alpha_1,\alpha_2,\alpha_3$ are inner angles of a Euclidean triangle with edge lengths $l_1,l_2,l_3$. By the cosine law, $\cos \alpha_i=\frac{-l_i^2+l_j^2+l_k^2}{2l_jl_k}$ for $\{i,j,k\}=\{1,2,3\}.$
Therefore $\alpha_1,\alpha_2,\alpha_3$ are functions in terms of $u_1,u_2,u_3.$

\begin{lemma}\label{thm:symmetry} The Jacobian matrix of functions $\alpha_1,\alpha_2,\alpha_3$ in terms of $u_1,u_2,u_3$ is symmetric.
\end{lemma}

\begin{proof} Due to the cosine law of a Euclidean triangles, $\alpha_1,\alpha_2,\alpha_3$ are functions of $l_1,l_2,l_3$. The following formula can be proved by direct calculation. For $\{i,j,k\}=\{1,2,3\}$,

$$\left(
\begin{array}{ccc}
d\alpha_1 \\
d\alpha_2 \\
d\alpha_3
\end{array}\right)
=\frac{-1}{\sin \alpha_il_jl_k} \left(
\begin{array}{ccc}
l_1&0&0 \\
0&l_2&0 \\
0&0&l_3
\end{array}
\right) \left(
\begin{array}{ccc}
-1&\cos \alpha_3&\cos \alpha_2 \\
\cos \alpha_3&-1&\cos \alpha_1 \\
\cos\alpha_2&\cos \alpha_1&-1
\end{array}
\right) \left(
\begin{array}{ccc}
dl_1 \\
dl_2 \\
dl_3
\end{array}\right).
$$

Next, by the construction $l_i^2=r_j^2+r_k^2+2I_{jk}r_jr_k,$ for $\{i,j,k\}=\{1,2,3\},$ where $I_{jk}$ is a constant. Differentiating the two sides of the equality, we obtained
$$dl_i=\frac{l_i^2+r_j^2-r_k^2}{2l_ir_j}dr_j+\frac{l_i^2+r_k^2-r_j^2}{2l_ir_k}dr_k.$$

Finally, $r_i=e^{u_i}$ implies $r_idu_i=dr_i$, for  $i=1,2,3.$

Combining the three relations, we have, for $\{i,j,k\}=\{1,2,3\}$,
\begin{align*}
\left(
\begin{array}{ccc}
d\alpha_1 \\
d\alpha_2 \\
d\alpha_3
\end{array}\right)
&=\frac{-1}{\sin \alpha_il_jl_k}  \left(
\begin{array}{ccc}
l_1&0&0 \\
0&l_2&0 \\
0&0&l_3
\end{array}
\right) \left(
\begin{array}{ccc}
-1&\cos \alpha_3&\cos \alpha_2 \\
\cos \alpha_3&-1&\cos \alpha_1 \\
\cos\alpha_2&\cos \alpha_1&-1
\end{array}
\right)\\
&\left(
\begin{array}{ccc}
0&\dfrac{l_1^2+r_2^2-r_3^2}{2l_1r_2}&\dfrac{l_1^2+r_3^2-r_2^2}{2l_1r_3} \\
\dfrac{l_2^2+r_1^2-r_3^2}{2l_2r_1}&0&\dfrac{l_2^2+r_3^2-r_1^2}{2l_2r_3} \\
\dfrac{l_3^2+r_1^2-r_2^2}{2l_3r_1}&\dfrac{l_3^2+r_2^2-r_1^2}{2l_3r_2}&0
\end{array}
\right)\left(
\begin{array}{ccc}
r_1&0&0 \\
0&r_2&0 \\
0&0&r_3
\end{array}
\right)\left(
\begin{array}{ccc}
du_1 \\
du_2 \\
du_3
\end{array}\right).
\end{align*}

We write the above formula as
$$
\left(
\begin{array}{ccc}
d\alpha_1 \\
d\alpha_2 \\
d\alpha_3
\end{array}\right)
=\frac{-1}{\sin \alpha_il_jl_k}N\left(
\begin{array}{ccc}
du_1 \\
du_2 \\
du_3
\end{array}\right).
$$

To show the Jacobian matrix of functions $\alpha_1,\alpha_2,\alpha_3$ in terms of $u_1,u_2,u_3$ is symmetric is equivalent to show the matrix $N$ is symmetric.

By the cosine law, we have
\begin{multline*}
4N=\left(
\begin{array}{ccc}
-2l_1^2 &l_1^2+l_2^2-l_3^2 &l_3^2+l_1^2-l_2^2 \\
l_1^2+l_2^2-l_3^2 &-2l_2^2 &l_2^2+l_3^2-l_1^2 \\
l_3^2+l_1^2-l_2^2 &l_2^2+l_3^2-l_1^2 &-2l_3^2
\end{array}
\right)
\left(
\begin{array}{ccc}
\frac1{l_1^2}&0&0 \\
0&\frac1{l_2^2}&0 \\
0&0&\frac1{l_3^2}
\end{array}
\right)\\
\left(
\begin{array}{ccc}
0 &l_1^2+r_2^2-r_3^2 &l_1^2+r_3^2-r_2^2 \\
l_2^2+r_1^2-r_3^2 &0 &l_2^2+r_3^2-r_1^2 \\
l_3^2+r_1^2-r_2^2 &l_3^2+r_2^2-r_1^2 &0
\end{array}
\right).
\end{multline*}

To simplify notations, let $a:=l_1^2,b:=l_2^2,c:=l_3^2, x:=\frac{r_2^2-r_3^2}a, y:=\frac{r_3^2-r_1^2}b, z:=\frac{r_1^2-r_2^2}c.$ Note that $ax+by+cz=0.$ Thus
\begin{align*}
4N=&\left(
\begin{array}{ccc}
-2a &a+b-c &c+a-b \\
a+b-c &-2b &b+c-a \\
c+a-b &b+c-a &-2c
\end{array}
\right)
\left(
\begin{array}{ccc}
0 &1+x &1-x \\
1-y &0 &1+y \\
1+z &1-z &0
\end{array}
\right)\\
&=\left(
\begin{smallmatrix}
2a-(a+b-c)y+(c+a-b)z &c-a-b-2ax-(c+a-b)z &b-a-c+2ax+(a+b-c)y \\
c-a-b+2by+(b+c-a)z &2b-(b+c-a)z+(a+b-c)x &a-b-c-2by-(a+b-c)x \\
b-c-a-2cz-(b+c-a)y &a-b-c+2cz+(c+a-b)x &2c-(c+a-b)x+(b+c-a)y
\end{smallmatrix}
\right).
\end{align*}
This matrix is symmetric due to the equality $ax+by+cz=0.$
\end{proof}

\begin{lemma}\label{thm:definite} The Jacobian matrix of functions $\alpha_1,\alpha_2,\alpha_3$ in terms of $u_1,u_2,u_3$ has one zero eigenvalue with associated eigenvector $(1,1,1)$ and two negative eigenvalues.
\end{lemma}

\begin{proof}
Using the notations in Lemma \ref{thm:symmetry}, we need to show that the matrix $4N$ has one zero eigenvalue and two positive eigenvalues. Since $4N$ is symmetric and the sum of entries in each row is zero, we can write $4N$ as
$$\left(
\begin{array}{ccc}
-B_2-B_3 &B_3 &B_2 \\
B_3 &-B_3-B_1 &B_1 \\
B_2 &B_1 &-B_1-B_2
\end{array}
\right),$$
where $B_1,B_2,B_3$ are functions of $a,b,c,x,y,z.$

The characteristic equation of $4N$ is $$\lambda(\lambda^2+2(B_1+B_2+B_3)\lambda+3(B_1B_2+B_2B_3+B_3B_1))=0.$$ There is a zero eigenvalue and the associated eigenvector is $(1,1,1)$. We claim that $B_1+B_2+B_3<0$ and $B_1B_2+B_2B_3+B_3B_1>0$. These two inequalities imply that the matrix $4N$ has two positive eigenvalues. In the following we verify the two inequalities.

First,
\begin{align*}
&  B_1+B_2+B_3<0\\
\Longleftrightarrow\  &  \mbox{trace of}\ 4N=-2(B_1+B_2+B_3)>0\\
\Longleftrightarrow\  &  a+b+c+(b-c)x+(c-a)y+(a-b)z>0\\
\Longleftrightarrow\  &  l_1^2+l_2^2+l_3^2+(l_2^2-l_3^2)\frac{r_2^2-r_3^2}{l_1^2}+
                                       (l_3^2-l_1^2)\frac{r_3^2-r_1^2}{l_2^2}+
                                       (l_1^2-l_2^2)\frac{r_1^2-r_2^2}{l_3^2}>0\\
\Longleftrightarrow\  &  l_1^2+l_2^2+l_3^2\\
                                     &+\frac{l_1^2l_2^2+l_1^2l_3^2-l_2^4-l_3^4}{l_2^2l_3^2}r_1^2
                                      +\frac{l_2^2l_1^2+l_2^2l_3^2-l_1^4-l_3^4}{l_1^2l_3^2}r_2^2
                                      +\frac{l_3^2l_1^2+l_3^2l_2^2-l_1^4-l_2^4}{l_1^2l_2^2}r_3^2>0\\
\Longleftrightarrow:\  & l_1^2+l_2^2+l_3^2 +\Xi_1r_1^2+\Xi_2r_2^2+\Xi_3r_3^2>0.
\end{align*}

For $\{i,j,k\}=\{1,2,3\}$, since $l_i^2=r_j^2+r_k^2+2I_{jk}r_jr_k$ and $I_{ij}\in[0,\infty),$ we have $l_i>r_j.$

Without loss of generality, we may assume $l_1\geq l_2\geq l_3.$ Therefore $\Xi_1\geq 0$ and $\Xi_3\leq 0.$

Case 1. If $\Xi_2\geq 0,$ then
\begin{align*}
&l_1^2+l_2^2+l_3^2 + \Xi_1r_1^2+\Xi_2r_2^2+\Xi_3r_3^2\\
&\geq l_1^2+l_2^2+l_3^2 + \Xi_3r_3^2\\
&>l_1^2+l_2^2+l_3^2 + \Xi_3l_2^2 \\
&=l_1^2+l_2^2+l_3^2 + \frac{l_3^2l_1^2+l_3^2l_2^2-l_1^4-l_2^4}{l_1^2l_2^2}l_2^2\\
&=\frac1{l_2^2}((l_1^2-l_2^2)l_2^2+2l_1^2l_3^2+l_2^2l_3^2)>0
\end{align*}
due to the assumption $l_1\geq l_2.$

Case 2. If $\Xi_2< 0,$ then
\begin{align*}
&l_1^2+l_2^2+l_3^2 + \Xi_1r_1^2+\Xi_2r_2^2+\Xi_3r_3^2\\
&>l_1^2+l_2^2+l_3^2 +\Xi_2r_2^2+\Xi_3r_3^2\\
&> l_1^2+l_2^2+l_3^2 +\Xi_2l_3^2+\Xi_3l_2^2\\
&= l_1^2+l_2^2+l_3^2 +\frac{l_2^2l_1^2+l_2^2l_3^2-l_1^4-l_3^4}{l_1^2l_3^2}l_3^2+\frac{l_3^2l_1^2+l_3^2l_2^2-l_1^4-l_2^4}{l_1^2l_2^2}l_2^2\\
&= \frac1{l_1^2}(2l_1^2l_2^2+2l_2^2l_3^2+2l_3^2l_1^2-l_1^4-l_2^4-l_3^4)\\
&= \frac1{l_1^2}(l_1+l_2+l_3)(l_1+l_2-l_3)(l_3+l_1-l_2)(l_2+l_3-l_1)>0
\end{align*}
due to the triangle inequality.

Thus $B_1+B_2+B_3<0$ holds.

To prove $B_1B_2+B_2B_3+B_3B_1>0$, we substitute back
$$B_1=4N_{23}=a-b-c-2by-(a+b-c)x,$$
$$B_2=4N_{31}=b-c-a-2cz-(b+c-a)y,$$
$$B_3=4N_{12}=c-a-b-2ax-(c+a-b)z $$ to obtain
\begin{align*}
&B_1B_2+B_2B_3+B_3B_1\\
=&\ 2ab+2bc+2ac-a^2-b^2-c^2\\
&+(b^2-c^2-3a^2+6ab+4ac)xy\\
&+(c^2-a^2-3b^2+6bc+4ab)yz\\
&+(a^2-b^2-3c^2+6ac+4bc)zx\\
&+2a(a+b-c)x^2+2b(b+c-a)y^2+2c(c+a-b)z^2\\
&+2(a+b+c)(ax+by+cz)\\
=&\ (2ab+2bc+2ac-a^2-b^2-c^2)(xy+yz+zx+1)\\
&+2(ax+by+cz)((a+b-c)x+(b+c-a)y+(c+a-b)z)\\
&+2(a+b+c)(ax+by+cz)\\
=&\ (2ab+2bc+2ac-a^2-b^2-c^2)(xy+yz+zx+1)
\end{align*}
due to the fact $ax+by+cz=0.$

Since
\begin{align*}
&2ab+2bc+2ac-a^2-b^2-c^2\\
=&\ 2l_1^2l_2^2+2l_2^2l_3^2+2l_3^2l_1^2-l_1^4-l_2^4-l_3^4\\
=&\ (l_1+l_2+l_3)(l_1+l_2-l_3)(l_3+l_1-l_2)(l_2+l_3-l_1)>0,
\end{align*}
we have
\begin{align}
&  B_1B_2+B_2B_3+B_3B_1>0\notag \\
\Longleftrightarrow\  & xy+yz+zx+1>0\notag \\
\Longleftrightarrow\  & (r_1^2-r_2^2)(r_3^2-r_1^2)l_1^2+(r_2^2-r_3^2)(r_1^2-r_2^2)l_2^2+(r_3^2-r_1^2)(r_2^2-r_3^2)l_3^2+l_1^2l_2^2l_3^2>0 \label{fml:e-1}
\end{align}

To verify inequality (\ref{fml:e-1}),  first note, since $l_i^2=r_j^2+r_k^2+2I_{jk}r_jr_k$ and $I_{jk}\in[0,\infty),$ we have $l_i^2\geq r_j^2+r_k^2$ for $\{i,j,k\}=\{1,2,3\}$.

Since the coefficient of $l_1^2$ in the left hand side of (\ref{fml:e-1}) satisfies
\begin{align*}
(r_1^2-r_2^2)(r_3^2-r_1^2)+l_2^2l_3^2&\geq (r_1^2-r_2^2)(r_3^2-r_1^2)+(r_1^2+r_3^2)(r_1^2+r_2^2)\\
&=2r_1^2r_2^2+2r_1^2r_3^2>0,
\end{align*}
by replacing $l_1^2$ by $r_2^2+r_3^2,$ we see that the left hand side of (\ref{fml:e-1}) is not less than
\begin{multline}\label{fml:e-2}
(r_1^2-r_2^2)(r_3^2-r_1^2)(r_2^2+r_3^2)+(r_2^2-r_3^2)(r_1^2-r_2^2)l_2^2+(r_3^2-r_1^2)(r_2^2-r_3^2)l_3^2\\+(r_2^2+r_3^2)l_2^2l_3^2.
\end{multline}

Since the coefficient of $l_2^2$ in (\ref{fml:e-2}) satisfies
\begin{align*}
(r_2^2-r_3^2)(r_1^2-r_2^2)+(r_2^2+r_3^2)l_3^2&\geq (r_2^2-r_3^2)(r_1^2-r_2^2)+(r_2^2+r_3^2)(r_1^2+r_2^2)\\
&=2r_1^2r_2^2+2r_2^2r_3^2>0,
\end{align*}
by replacing $l_2^2$ by $r_1^2+r_3^2,$ we see that (\ref{fml:e-2}) is not less than
\begin{multline}\label{fml:e-3}
(r_1^2-r_2^2)(r_3^2-r_1^2)(r_2^2+r_3^2)+(r_2^2-r_3^2)(r_1^2-r_2^2)(r_1^2+r_3^2)\\+(r_3^2-r_1^2)(r_2^2-r_3^2)l_3^2+(r_2^2+r_3^2)(r_1^2+r_3^2)l_3^2.
\end{multline}

Since the coefficient of $l_3^2$ in (\ref{fml:3}) is positive, by replacing $l_3^2$ by $r_1^2+r_2^2,$ we see that (\ref{fml:e-3}) is not less than
\begin{multline}\label{fml:e-4}
(r_1^2-r_2^2)(r_3^2-r_1^2)(r_2^2+r_3^2)+(r_2^2-r_3^2)(r_1^2-r_2^2)(r_1^2+r_3^2)\\
+(r_3^2-r_1^2)(r_2^2-r_3^2)(r_1^2+r_2^2)+(r_2^2+r_3^2)(r_1^2+r_3^2)(r_1^2+r_2^2).
\end{multline}
Since (\ref{fml:e-4}) is equal to $8r_1^2r_2^2r_3^2$, we see that  $B_1B_2+B_2B_3+B_3B_1>0.$

\end{proof}

For the special case of intersection angle circle packing, Lemma \ref{thm:definite} was proved in Chow-Luo \cite{cl}(Lemma 3.1). In that case, due to Thurston \cite{th}, the monotonicity of angles holds i.e., $\frac{\partial \alpha_i}{\partial r_i}<0, \frac{\partial \alpha_i}{\partial r_j}>0$ for $j\neq i.$ Combining with linear algebra, this property implies Lemma \ref{thm:definite}.
But in the general case of inversive distance circle packing, the geometric picture does not hold. For example, when $l_1=2,l_2=2,l_3=3,r_1=r_2=r_3=1,$ the matrix $N$ is
$$\left(
\begin{array}{ccc}
2 &\frac14 &-\frac94 \\
\frac14 &2 &-\frac94 \\
-\frac94 &-\frac94 &\frac92
\end{array}
\right)$$ with eigenvalue $0, \frac{27}4,\frac74.$

Since, by Corollary \ref{thm:e-space}, the space $\mathcal{U}_E^{123}$ of vectors $(u_1,u_2,u_3)$ is connected and simply connected, Lemma \ref{thm:symmetry} and Lemma \ref{thm:definite} imply

\begin{corollary}\label{thm:functional}
The differential 1-form $\sum_{i=1}^3\alpha_idu_i$ is closed. For any $c\in \mathcal{U}_E^{123},$
the integration $w(u_1,u_2,u_3)=\int^{(u_1,u_2,u_3)}_c\sum_{i=1}^3\alpha_idu_i$
is a concave function on $\mathcal{U}_E^{123}$ and satisfying, for $i=1,2,3,$
\begin{align}\label{fml:gradient}
\frac{\partial w}{\partial u_i} =\alpha_i.
\end{align}
\end{corollary}

\begin{proof}[Proof of Theorem \ref{thm:main} in Euclidean geometry] Let's prove the local rigidity of inversive distance circle packing in Euclidean geometry. For a closed triangulated weighted surface
$(\Sigma,T, I)$ with $I_{ij}\in [0,\infty)$ for each $ij\in E,$ let $\mathcal{U}_E$ be the open subset of $\mathbb{R}^{|V|}$ formed by the vectors $u=(u_1,u_2,..,u_{|V|})$ satisfying $(u_i,u_j,u_k)\in \mathcal{U}_E^{ijk}$ whenever $\triangle ijk \in F.$

For a radius vector $r=(r_1,r_2,..,r_{|V|})$, we may rescale it to get $cr$ while the cone angles remain the same. Under the rescaling, $u$ changes to $u+\ln(c)(1,1,...,1).$ Therefore we only consider the intersection $\mathcal{U}_E\cap P$ where $P$ is the hyperplane defined by $\sum_{i=1}^{|V|}u_i=0.$

Let $a_i$ be the cone angle at vertex $i$ under the Euclidean cone metric. Then $\sum_{i=1}^{|V|}a_i=$ the sum of all inner angles $=\pi|F|.$ Therefore the space of possible cone angles is contained in the intersection $\mathbb{R}^{|V|}_{>0}\cap Q$ where $Q$ is the hyperplane defined by $\sum_{i=1}^{|V|}a_i=\pi|F|.$ We have a map $\xi:\mathcal{U}_E\cap P\to \mathbb{R}^{|V|}_{>0}\cap Q.$

By Corollary \ref{thm:functional}, for each triangle $\triangle ijk\in F$, there is a function $w(u_i,u_j,u_k).$ Define a function $W: \mathcal{U}_E \to\mathbb{R}$ by
$$W(u_1,u_2,...,u_{|V|})=\sum_{\triangle ijk \in F} w(u_i, u_j, u_k)$$ where the sum is
over all triangles in $F$.

By (\ref{fml:gradient}), $\frac{\partial W}{\partial u_i}$ equals the sum of inner angles having vertex $i$, i.e., the cone angle at vertex $i$. Therefore, when restricted on $\mathcal{U}_E\cap P$, the gradient $\nabla W$ is the map $\xi:\mathcal{U}_E\cap P\to \mathbb{R}^{|V|}_{>0}\cap Q.$

We claim that the restriction of $W$ on $\mathcal{U}_E\cap P$ is strictly concave. By the definition of $W,$ we get the relation between the Hessian matrixes $$H(W)=\sum_{\triangle ijk \in F} H(w(u_i, u_j, u_k)),$$ where we think $w(u_i, u_j, u_k)$ as a function on $\mathcal{U}_E.$

For any $x \in \mathcal{U}_E\cap P,$
\begin{align*}
x H(W)x^T &=\sum_{\triangle ijk \in F} x H(w(u_i, u_j, u_k))x^T\\
&= \sum_{\triangle ijk \in F} (x_i,x_j,x_k) H(w(u_i, u_j, u_k))(x_i,x_j,x_k)^T \leq 0
\end{align*}
since $H(w(u_i, u_j, u_k))$ is negative semi-definite for each $\triangle ijk \in F.$

Thus $x H(W)x^T=0$ implies $(x_i,x_j,x_k) H(w(u_i, u_j, u_k))(x_i,x_j,x_k)^T=0$ for each $\triangle ijk.$ By Lemma \ref{thm:definite}, $(x_i,x_j,x_k)$ is an eigenvector of zero eigenvalue. Therefore $(x_i,x_j,x_k)=c_{ijk}(1,1,1)$ for some constant $c_{ijk}$. If two triangles $\triangle ijk$ and $\triangle i'j'k'$ share a vertex, then $c_{ijk}=c_{i'j'k'}.$ This implies that $c_{ijk}$ must be a constant $c$ independent of triangles $\triangle ijk.$ Therefore $x=c(1,1,...,1).$ Since $x\in P$, i.e, $\sum x_i=0,$ we have $c=0.$ Thus $x=(0,0,...,0)$. This shows that the Hessian $H(W)$ is negative definite when restricted on $\mathcal{U}_E\cap P$.

The local rigidity of inversive distance circle packing in Euclidean geometry follows from the following lemma in analysis

\begin{lemma}\label{thm:convex} If $X$ is an open set in $\mathbb{R}^n$
and the Hessian matrix of $f$ is positive definite for all $x \in X$, then the
gradient $\nabla f:X\to\mathbb{R}^n$ is a local diffeomorphism.
\end{lemma}

\end{proof}

\section{Circle packing in hyperbolic geometry}

\subsection{Space of radius vectors}

As stated in Introduction, to obtain a hyperbolic cone metric, we need conditions on a radius vector.
In this section we investigate the space of radius vectors for one triangle in hyperbolic geometry.
Consider a topological triangle
$\triangle ijk\in F$. Assign a number $l_{ij}$ to the edge $ij\in E,$ where $l_{ij}$ is the distance of the centers of two circles with radii $r_i,r_j$ and inversive distance $I_{ij}.$ Therefore $l_{ij}$ is a function of $r_i, r_j.$ In fact, in hyperbolic geometry, $l_{ij}$ is determined by
$\cosh l_{ij}=\cosh r_i\cosh r_j+I_{ij}\sinh r_i\sinh r_j.$ Similarly, we obtain $l_{jk},l_{ki}.$

\begin{lemma} The set $\mathcal{R}_H^{ijk}:=\{(r_i, r_j,r_k)\in \mathbb{R}^3_{>0}\ |\ l_{ij},l_{jk},l_{ki}$ satisfy triangle inequality\} is a connected and simply connected open subset in $\mathbb{R}^3_{>0}.$
\end{lemma}

\begin{proof} 
Let $\rho_H:\mathbb{R}^3_{> 0}\to \mathbb{R}^3_{> 0}$ be the map sending $(r_i,r_j,r_k)$ to $(l_{ij},l_{jk},l_{ki}).$ We claim that $\rho_H$ is a smooth embedding. First the determinant of the Jacobian matrix of $\rho_H$ is nonzero. Next, we will show $\rho_H$ is one-to-one.

That means, for fixed numbers $I_{ij},I_{jk},I_{ki}\in [0,\infty)$ and a fixed hyperbolic triangle, there is at most one configuration of three circles centered at the vertexes of the triangle and having the prescribed inversive distance $I_{ij},I_{jk},I_{ki}$.

The statement is true in Euclidean geometry as proved in the proof of Lemma \ref{thm:e-r-space}. We use the result in Euclidean geometry to prove the similar result in hyperbolic geometry. Let's assume there are two configurations of three circles in hyperbolic plane with the fixed triangle and inversive distance. As in Figure \ref{fig:projection}, the stereographic projection sends the two configurations to Euclidean plane $\mathbb{C}.$ Since the stereographic projection preserve the inversive distance, there will be two different configuration of three circles with the same Euclidean triangle and the same inversive distance. This is impossible.

To determine the image $\rho_H(\mathbb{R}^3_{> 0}),$ we consider the image of the boundary of $\mathbb{R}^3_{\geq 0}.$ If $r_i=0,$ then $r_k=l_{ki}, r_j=l_{ij}.$ Thus $$\cosh l_{jk}=\cosh r_j\cosh r_k+I_{jk}\sinh r_j\sinh r_k=\cosh l_{ij}\cosh l_{ki}+I_{jk}\sinh l_{ij}\sinh l_{ki}.$$ On the other hand, if $r_i>0,$ then $r_k<l_{ki}, r_j<l_{ij}.$ We have $$\cosh l_{jk}<\cosh l_{ij}\cosh l_{ki}+I_{jk}\sinh l_{ij}\sinh l_{ki}.$$ Therefore the image $\rho_H(\mathbb{R}^3_{> 0})$ is the subset of vectors $(l_{ij},l_{jk},l_{ki})$ satisfying
$$\cosh l_{jk}<\cosh l_{ij}\cosh l_{ki}+I_{jk}\sinh l_{ij}\sinh l_{ki},$$
$$\cosh l_{ki}<\cosh l_{jk}\cosh l_{ij}+I_{ki}\sinh l_{jk}\sinh l_{ij},$$
$$\cosh l_{ij}<\cosh l_{ki}\cosh l_{jk}+I_{ij}\sinh l_{ki}\sinh l_{jk}.$$
In fact, $\rho_H(\mathbb{R}^3_{> 0})$ is a cone bounded by three surfaces. In $\mathbb{R}^3_{\geq 0}$, the intersection of any of the two surfaces is a straight line. For example, if 
$$\cosh l_{jk}=\cosh l_{ij}\cosh l_{ki}+I_{jk}\sinh l_{ij}\sinh l_{ki}$$ and 
$$\cosh l_{ki}=\cosh l_{jk}\cosh l_{ij}+I_{ki}\sinh l_{jk}\sinh l_{ij},$$ then the sum of the two equation gives
\begin{align*}
0&=(\cosh l_{ki}+\cosh l_{jk})(\cosh l_{ij}-1)+(I_{jk}\sinh l_{ki}+I_{ki}\sinh l_{jk})\sinh l_{ij}\\
&=(\cosh l_{ki}+\cosh l_{jk})2\sinh^2\frac{l_{ij}}{2}+(I_{jk}\sinh l_{ki}+I_{ki}\sinh l_{jk})2\cosh\frac{l_{ij}}{2}\sinh\frac{l_{ij}}{2}.
\end{align*}
The only possibility is $\sinh\frac{l_{ij}}{2}=0$. Therefore $l_{ij}=0$ and $l_{jk}=l_{ki}.$ 

Hence the three lines of intersection are 
$$l_{ij}=0 \ \mbox{and}\ l_{jk}=l_{ki},$$
$$l_{jk}=0 \ \mbox{and}\ l_{ki}=l_{ij},$$
$$l_{ki}=0 \ \mbox{and}\ l_{ij}=l_{jk}.$$

Recall that $\mathcal{L}$ is the subset of $\mathbb{R}^3_{> 0}$ formed by vectors $(l_{ij},l_{jk},l_{ki})$ satisfying the triangle inequality. Therefore $\mathcal{L}$ is a cone bounded by three planes. The three lines of intersection of the three planes are the same as that of the boundary surfaces of $\rho_H(\mathbb{R}^3_{> 0}).$ 

We claim that $\rho_H(\mathbb{R}^3_{> 0})\cap \mathcal{L}$ is a cone.

In fact, if $I_{ij}\in [0,1],$ then, in $\mathbb{R}^3_{> 0}$, we have
\begin{multline*}
\{(l_{ij},l_{jk},l_{ki})\ |\ \cosh l_{ij}<\cosh l_{ki}\cosh l_{jk}+I_{ij}\sinh l_{ki}\sinh l_{jk}\}\\
\subseteq\{(l_{ij},l_{jk},l_{ki})\ |\ l_{ij}< l_{ki}+l_{jk}\}.
\end{multline*}
If $I_{ij}\in (1,\infty),$ then, in $\mathbb{R}^3_{> 0}$, we have
\begin{multline*}
\{(l_{ij},l_{jk},l_{ki})\ |\ \cosh l_{ij}<\cosh l_{ki}\cosh l_{jk}+I_{ij}\sinh l_{ki}\sinh l_{jk}\}\\
\supset\{(l_{ij},l_{jk},l_{ki})\ |\ l_{ij}< l_{ki}+l_{jk}\}.
\end{multline*}

Applying the same argument in Euclidean geometry, we see $\rho_H(\mathbb{R}^3_{> 0})\cap \mathcal{L}$ is a cone. Therefore 
$\mathcal{R}_H^{ijk}=\rho_H^{-1}(\rho_H(\mathbb{R}^3_{> 0})\cap \mathcal{L})$ is connected and simply connected.

\end{proof}

In the construction of the energy function, we will use a natural variable $u=(u_1,u_2,...,u_{|V|})$ where $u_i= \ln \tanh\frac{r_i}2.$ Let $\tau_H:\mathbb{R}^3 \to \mathbb{R}^3_{>0}$ be the map sending $(u_i,u_j,u_k)$ to $(r_i,r_j,r_k).$ Since $\tau_H$ is a diffeomorphism, we have

\begin{corollary}\label{thm:h-space}
The set $\mathcal{U}_H^{ijk}:=\tau_H^{-1}(\mathcal{R}_H^{ijk})$ is a connected and simply connected open subset in $\mathbb{R}^3_{>0}.$
\end{corollary}

To prove Theorem \ref{thm:main}, we construct a strictly concave energy function on $\mathcal{U}_H^{ijk}$ for each triangle $\triangle ijk\in F.$
As in the case of Euclidean geometry, the space $\mathcal{U}_H^{ijk}$ is not always convex (for example, when $I_{ij}=I_{jk}=I_{ki}=2$), we can not obtain the global rigidity. Therefore we are not able to determine the space of cone angles in the case of inversive distance circle packing.

\subsection{Energy function}

To prove Theorem \ref{thm:main}, we will construct an energy function on the space of radius vectors. First, we construct an energy function on the space for one triangle. To simplify notations, consider the triangle $\triangle 123\in F$ and the map
$$\begin{array}{ccccccc}
\mathcal{U}_H^{123} & \to & \mathcal{R}_H^{123}& \to &  \mathcal{L}& \to &  \mathbb{R}^3_{>0}  \\
(u_1,u_2,u_3) & \mapsto &  (r_1,r_2,r_3) & \mapsto &  (l_1,l_2,l_3) & \mapsto &  (\alpha_1,\alpha_2,\alpha_3)
\end{array}$$
where $r_i=\frac{1+e^{u_i}}{1-e^{u_i}}$ for $i=1,2,3.$ For three fixed positive numbers $I_{12},I_{23},I_{31},$ we have $\cosh l_i=\cosh r_j\cosh r_k+I_{ij}\sinh r_j\sinh r_k$ for $\{i,j,k\}=\{1,2,3\}.$ For simplicity of notations, we set $(l_1,l_2,l_3)=(l_{23},l_{31},l_{12}).$ And $\alpha_1,\alpha_2,\alpha_3$ are inner angles of a hyperbolic triangle with edge lengths $l_1,l_2,l_3$. By the cosine law, $$\cos \alpha_i=\frac{-\cosh l_i+\cosh l_j\cosh l_k}{\sinh l_j\sinh l_k}$$ for $\{i,j,k\}=\{1,2,3\}.$
Therefore $\alpha_1,\alpha_2,\alpha_3$ are functions in terms of $u_1,u_2,u_3.$

\begin{lemma}\label{thm:h-symmetry} The Jacobian matrix of functions $\alpha_1,\alpha_2,\alpha_3$ in terms of $u_1,u_2,u_3$ is symmetric.
\end{lemma}

\begin{proof} Due to the cosine law of a hyperbolic triangles, $\alpha_1,\alpha_2,\alpha_3$ are functions of $l_1,l_2,l_3$. The following formula can be proved by direct calculation. For $\{i,j,k\}=\{1,2,3\}$,

\begin{multline*}
\left(
\begin{array}{ccc}
d\alpha_1 \\
d\alpha_2 \\
d\alpha_3
\end{array}\right)
=\frac{-1}{\sin \alpha_i\sinh l_j\sinh l_k} \left(
\begin{array}{ccc}
\sinh l_1&0&0 \\
0&\sinh l_2&0 \\
0&0&\sinh l_3
\end{array}
\right)\\
 \left(
\begin{array}{ccc}
-1&\cos \alpha_3&\cos \alpha_2 \\
\cos \alpha_3&-1&\cos \alpha_1 \\
\cos\alpha_2&\cos \alpha_1&-1
\end{array}
\right) \left(
\begin{array}{ccc}
dl_1 \\
dl_2 \\
dl_3
\end{array}\right).
\end{multline*}

Next, by the construction $\cosh l_i=\cosh r_j\cosh r_k+I_{jk}\sinh r_j\sinh r_k,$ for $\{i,j,k\}=\{1,2,3\},$ where $I_{ij}$ is a constant. Differentiating the two sides of the equality, we obtained
$$dl_i=\frac{-\cosh r_k+\cosh l_i\cosh r_j}{\sinh l_i\sinh r_j}dr_j+
\frac{-\cosh r_j+\cosh l_i\cosh r_k}{\sinh l_i\sinh r_k}dr_k.$$

Finally, $r_i=\ln\frac{1+e^{u_i}}{1-e^{u_i}}$ implies $\sinh r_idu_i=dr_i, i=1,2,3.$

Combining the three relations, we have, for $\{i,j,k\}=\{1,2,3\}$,

\begin{multline*}
\left(
\begin{array}{ccc}
d\alpha_1 \\
d\alpha_2 \\
d\alpha_3
\end{array}\right)
=\frac{-1}{\sin \alpha_i\sinh l_j\sinh l_k}\\
\left(
\begin{array}{ccc}
\sinh l_1&0&0 \\
0&\sinh l_2&0 \\
0&0&\sinh l_3
\end{array}
\right)
\left(
\begin{array}{ccc}
-1&\cos \alpha_3&\cos \alpha_2 \\
\cos \alpha_3&-1&\cos \alpha_1 \\
\cos\alpha_2&\cos \alpha_1&-1
\end{array}
\right)\\
\left(
\begin{array}{ccc}
0&\dfrac{-\cosh r_3+\cosh l_1\cosh r_2}{\sinh l_1\sinh r_2}&\dfrac{-\cosh r_2+\cosh l_1\cosh r_3}{\sinh l_1\sinh r_3} \\
\dfrac{-\cosh r_3+\cosh l_2\cosh r_1}{\sinh l_2\sinh r_1}&0&\dfrac{-\cosh r_1+\cosh l_2\cosh r_3}{\sinh l_2\sinh r_3} \\
\dfrac{-\cosh r_2+\cosh l_3\cosh r_1}{\sinh l_3\sinh r_1}&\dfrac{-\cosh r_1+\cosh l_3\cosh r_2}{\sinh l_3\sinh r_2}&0
\end{array}
\right)\\
\left(
\begin{array}{ccc}
\sinh r_1&0&0 \\
0&\sinh r_2&0 \\
0&0&\sinh r_3
\end{array}
\right)\left(
\begin{array}{ccc}
du_1 \\
du_2 \\
du_3
\end{array}\right).
\end{multline*}

We write the above formula as
$$
\left(
\begin{array}{ccc}
d\alpha_1 \\
d\alpha_2 \\
d\alpha_3
\end{array}\right)
=\frac{-1}{\sin \alpha_i\sinh l_j\sinh l_k}M\left(
\begin{array}{ccc}
du_1 \\
du_2 \\
du_3
\end{array}\right),
$$
where $M$ is a product of four matrixes.

To show the Jacobian matrix of functions $\alpha_1,\alpha_2,\alpha_3$ in terms of $u_1,u_2,u_3$ is symmetric is equivalent to show the matrix $M$ is symmetric.

To simplify notations, let $a:=\cosh l_1,b:=\cosh l_2,c:=\cosh l_3, x:=\cosh r_1, y:=\cosh r_2, z:=\cosh r_3.$
By the cosine law, we have
\begin{align*}
M=&\left(
\begin{array}{ccc}
1-a^2 &ab-c &ca-b \\
ab-c &1-b^2 &bc-a \\
ca-b &bc-a &1-c^2
\end{array}
\right)\left(
\begin{array}{ccc}
\frac1{a^2-1}&0&0 \\
0&\frac1{b^2-1}&0 \\
0&0&\frac1{c^2-1}
\end{array}
\right)
\left(
\begin{array}{ccc}
0 &ay-z &az-y \\
bx-z &0 &bz-x \\
cx-y &cy-x &0
\end{array}
\right).
\end{align*}

We check that $M$ is symmetric, for example, by showing $$M_{12}=M_{21}=z-\frac{ac-b}{c^2-1}x-\frac{bc-a}{c^2-1}y.$$
\end{proof}

\begin{lemma}\label{thm:h-definite} The Jacobian matrix of functions $\alpha_1,\alpha_2,\alpha_3$ in terms of $u_1,u_2,u_3$ is negative definite.
\end{lemma}

\begin{proof} Using the notations of Lemma \ref{thm:h-symmetry}, we need to verify that the matrix $M$ is positive definite. The determinant of $M$ is positive since the determinant of each factor is positive. Note that the second matrix is a Gram matrix up to the negative sign. Since, by Corollary \ref{thm:h-space}, the set $\mathcal{U}_H^{123}$ is connected, to show $M$ is positive
definite, it is enough to check $M$ is positive definite at one
vector in $\mathcal{U}_H^{123}$.

It is equivalent to pick up a vector $(r_1,r_2,r_3)$ such that the resulting $l_1,l_2,l_3$ satisfy the triangle inequality. We pick up $(r_1, r_2, r_3)=(s,s,s)$. In the following, we show that such kind of number $s$ exists.

For fixed $I_{12}, I_{23}, I_{31}\in[0,\infty),$ we have $\cosh l_i=\cosh r_j\cosh r_k+I_{jk}\sinh r_j\sinh r_k=\cosh^2 s+I_{jk}\sinh^2 s \geq \cosh^2 s,$ for $\{i,j,k\}=\{1,2,3\}.$ Therefore $\sinh l_i\geq \sqrt{\cosh^4 s-1}$ for $i=1,2,3.$

That $l_1,l_2,l_3$ satisfy the triangle inequality is equivalent to

\begin{align*}
&\cosh(l_i+l_j)>\cosh l_k\\
\Longleftrightarrow& \cosh l_i\cosh l_j+\sinh l_i\sinh l_j>\cosh l_k=\cosh^2 s+I_{ij}\sinh^2 s\\
\Longleftarrow& \cosh^4 s+\cosh^4 s-1> \cosh^2 s+I_{ij}\sinh^2 s.
\end{align*}
Let $x:=\cosh^2 s.$ The above inequality is equivalent to $(2x-I_{ij}+1)(x-1)>0.$ Since $x>1,$ we need $x>\frac{I_{ij}-1}2.$

Therefore $l_1,l_2,l_3$ as functions of $(s,s,s)$ satisfy the triangle inequality if $\cosh^2 s>\max\{\frac{I_{12}-1}2, \frac{I_{23}-1}2, \frac{I_{31}-1}2\}.$

Now it is enough to check that $M$ is positive definite at such a vector $(s,s,s)$ when $s$ is sufficiently large. In this case,
\begin{multline*}
M=\left(
\begin{array}{ccc}
\sinh l_1&0&0 \\
0&\sinh l_2&0 \\
0&0&\sinh l_3
\end{array}
\right)
\left(
\begin{array}{ccc}
-1&\cos \alpha_3&\cos \alpha_2 \\
\cos \alpha_3&-1&\cos \alpha_1 \\
\cos\alpha_2&\cos \alpha_1&-1
\end{array}
\right)\\
\left(
\begin{array}{ccc}
0&\dfrac{-1+\cosh l_1}{\sinh l_1}&\dfrac{-1+\cosh l_1}{\sinh l_1} \\
\dfrac{-1+\cosh l_2}{\sinh l_2}&0&\dfrac{-1+\cosh l_2}{\sinh l_2} \\
\dfrac{-1+\cosh l_3}{\sinh l_3}&\dfrac{-1+\cosh l_3}{\sinh l_3}&0
\end{array}
\right)\cosh s:=M_1\cosh s.
\end{multline*}

$M$ is positive definite is equivalent to that $M_1$ is positive definite. Since $M_1$ only involves $l_1,l_2,l_3$, that $M_1$ is positive definite is a property of a hyperbolic triangle. In fat, that $M_1$ is positive definite was proved in the paper Guo-Luo \cite{gl} (Lemma 4.4), by showing that its leading principal $1\times 1$ and $2\times 2$ minor are positive. For completeness, we include a proof here.

That the leading principal $1\times 1$ minor is positive is equivalent to
$$\cos \alpha_3\frac{-1+\cosh l_2}{\sinh l_2}+\cos \alpha_2\frac{-1+\cosh l_3}{\sinh l_3}>0.$$
Replacing $\alpha'$s by $l'$s by the cosine law, we obtain
\begin{equation}\label{fml:3}
\cosh l_1>\frac{\cosh^2 l_2+\cosh^2 l_3+\cosh l_2+\cosh
l_3}{2\cosh l_2\cosh l_3+\cosh l_2+\cosh l_3}.
\end{equation}

Since $l_1>|l_2-l_3|,$ therefore $\cosh l_1>\cosh (l_2-l_3).$ To
show (\ref{fml:3}) holds, it is enough to check
\begin{equation}\label{fml:4}
\cosh l_2\cosh l_3-\sinh l_2\sinh l_3\geq \frac{\cosh^2
l_2+\cosh^2 l_3+\cosh l_2+\cosh l_3}{2\cosh l_2\cosh l_3+\cosh
l_2+\cosh l_3}.
\end{equation}

We simplify the notations by introducing $a:=\cosh l_2>1, b:=\cosh
l_3>1.$ Then (\ref{fml:4}) is rewritten as
\begin{equation}\label{fml:5}
ab-\sqrt{(a^2-1)(b^2-1)}\geq\frac{a^2+b^2+a+b}{2ab+a+b}.
\end{equation}
(\ref{fml:5}) is equivalent to
\begin{equation}\label{fml:6}
ab-\frac{a^2+b^2+a+b}{2ab+a+b}\geq\sqrt{(a^2-1)(b^2-1)}.
\end{equation}

Since $a>1,b>1,$ the left hand side of (\ref{fml:6}) is positive.
To show (\ref{fml:6}) holds, we square the two sides and simplify.
We have
\begin{align*}
&(ab^4+a^4b-a^3b^2-a^2b^3)+(a^4+b^4-2a^2b^2)+(a^3+b^3-ab^2-a^2b)\geq0\\
\Longleftrightarrow& (ab+1)(a^3+b^3-a^2b-ab^2)+(a^2-b^2)^2\geq0\\
\Longleftrightarrow& (ab+1)(a+b)(a-b)^2+(a^2-b^2)^2\geq0.
\end{align*}
This shows that the leading principal $1\times 1$ minor is
positive.

To show the leading principal $2\times 2$ minor of $M_1$ is
positive, for simplifying notations, let $t_i:=\frac{-1+\cosh l_i}{\sinh l_i}$ . Up to a diagonal matrix, $M_1$ is
$$
\left(
\begin{array}{ccc}
-1&\cos \alpha_3&\cos \alpha_2 \\
\cos \alpha_3&-1&\cos \alpha_1 \\
\cos\alpha_2&\cos \alpha_1&-1
\end{array}
\right)
\left(
\begin{array}{ccc}
0&t_1&t_1 \\
t_2&0&t_2 \\
t_3&t_3&0
\end{array}
\right).
$$

That the leading principal $2\times 2$ minor of $M_1$ is
positive is equivalent to
\begin{align*}
&(t_2\cos \alpha_3+t_3\cos \alpha_2)(t_1\cos \alpha_3+t_3\cos \alpha_1)-(-t_2+t_3\cos \alpha_1)(-t_1+t_3\cos \alpha_2)>0\\
\Longleftrightarrow\  & (\cos^2 \alpha_3-1)t_1t_2+(\cos \alpha_1\cos \alpha_3+\cos \alpha_2)t_2t_3+(\cos \alpha_2\cos \alpha_3+\cos \alpha_1)t_1t_3>0\\
\Longleftrightarrow\  &  -\sin^2 \alpha_3 t_1t_2+\cosh l_2 \sin \alpha_1\sin \alpha_3t_2t_3+\cosh l_1 \sin \alpha_2\sin \alpha_3t_1t_3>0\\
\Longleftrightarrow\  & \cosh l_2 \sinh l_1t_2t_3+\cosh l_1 \sinh l_2t_1t_3>\sinh l_3 t_1t_2\\
\Longleftrightarrow\  & \frac{\cosh l_2 \sinh l_1}{t_1}+\frac{\cosh l_1 \sinh l_2}{t_2}>\frac{\sinh l_3}{t_3}\\
\Longleftrightarrow\  & \cosh l_2(1+\cosh
l_1)+ \cosh l_1(1+\cosh l_2)>1+\cosh l_3 \\
\Longleftrightarrow\  & (\cosh l_1+\cosh l_2-1)+(2\cosh l_1\cosh
l_2-\cosh l_3)>0.
\end{align*}
The cosine law and the sine law of a hyperbolic triangle are used in above calculation.

The last inequality is true since $\cosh l_1+\cosh l_2>1$ and $2\cosh
l_1\cosh l_2>\cosh l_1\cosh l_2+\sinh l_1\sinh
l_2=\cosh(l_1+l_2)>\cosh l_3.$

\end{proof}

For the special case of intersection angle circle packing, Lemma \ref{thm:h-definite} was proved in Chow-Luo \cite{cl}(Lemma 3.1) by applying the geometric picture of monotonicity of angles (an observation due to Thurston \cite{th}).
But in the general case of inversive distance circle packing, the geometric picture does not hold. For example, when $l_1=2,l_2=2,l_3=3,r_1=r_2=r_3=1,$ the matrix $M$ is
$$\left(
\begin{array}{ccc}
6.08 &0.49 &-2.94 \\
0.49 &6.08 &-2.94 \\
-2.94 &-2.94 &22.11
\end{array}
\right)$$ with eigenvalue $23.15, 5.59, 5.53.$

Since, by Corollary \ref{thm:h-space}, the space $\mathcal{U}_H^{123}$ of vectors $(u_1,u_2,u_3)$ is connected and simply connected, Lemma \ref{thm:h-symmetry} and Lemma \ref{thm:h-definite} imply

\begin{corollary}\label{thm:h-functional}
The differential 1-form $\sum_{i=1}^3\alpha_idu_i$ is closed. For any $c\in \mathcal{U}_H^{123},$
the integration $w(u_1,u_2,u_3)=\int^{(u_1,u_2,u_3)}_c\sum_{i=1}^3\alpha_idu_i$
is a strictly concave function on $\mathcal{U}_H^{123}$ and satisfying, for $i=1,2,3,$
\begin{align}\label{fml:h-gradient}
\frac{\partial w}{\partial u_i} =\alpha_i.
\end{align}
\end{corollary}

\begin{proof}[Proof of Theorem \ref{thm:main} in hyperbolic geometry] Let's prove the local rigidity of inversive distance circle packing in hyperbolic geometry. For a closed triangulated weighted surface
$(\Sigma,T, I)$ with $I_{ij}\in [0,\infty)$ for each $ij\in E,$ let $\mathcal{U}_H$ be the open subset of $\mathbb{R}^{|V|}$ formed by the vectors $u=(u_1,u_2,..,u_{|V|})$ satisfying $(u_i,u_j,u_k)\in \mathcal{U}_H^{ijk}$ whenever $\triangle ijk \in F.$

By Corollary \ref{thm:h-functional}, for each triangle $\triangle ijk\in F$, there is a function $w(u_i,u_j,u_k).$ Define a function $W: \mathcal{U}_H \to\mathbb{R}$ by
$$W(u_1,u_2,...,u_{|V|})=\sum_{\triangle ijk \in F} w(u_i, u_j, u_k)$$ where the sum is
over all triangles in $F$. By
Corollary \ref{thm:h-functional}, $W$ is strictly concave on
$\mathcal{U}_H$ so that
$\frac{\partial W}{\partial u_i}$ equals the sum of inner angles having vertex $i$, i.e., the cone angle at vertex $i$. That means the gradient of $W$ is exactly the map sending a vector $u$ to its cone angles.

The local rigidity of inversive distance circle packing in hyperbolic geometry holds due to Lemma \ref{thm:convex}.

\end{proof}

\section*{Acknowledgment} 

The author would like to thank Feng Luo and Albert Madern for encouragement and helpful comments.

\bibliographystyle{amsplain}

\end{document}